\documentclass[10pt,oneside,english]{amsart}

\usepackage{amsthm}
\usepackage{amssymb}

\makeatletter
\numberwithin{equation}{section} 
\numberwithin{figure}{section} 
\theoremstyle{plain}
\theoremstyle{plain}
\newtheorem{thm}{Theorem}
  \theoremstyle{definition}
  \newtheorem{defn}[thm]{Definition}
  \theoremstyle{plain}
  \newtheorem{lem}[thm]{Lemma}
  \theoremstyle{remark}
  \newtheorem{rem}[thm]{Remark}
  \theoremstyle{plain}
  \newtheorem{cor}[thm]{Corollary}

\makeatother

\usepackage{babel}

\begin{document}
\textheight=8.5in \textwidth=6.0in \topmargin=-0.3in

\title[traveling wave solution]{Traveling waves and their Stability for a Public Goods Game Model }

\maketitle
\global\long\def\theequation{\arabic{section}.\arabic{equation}}


\begin{center}
XIAOJIE HOU, \textsc{WEI FENG}%
\footnote{Contacting author

\medskip{}


} 
\par\end{center}

{\footnotesize \centerline{Department of Mathematics and Statistics,}
\centerline{ University of North Carolina at Wilmington,} \centerline{
Wilmington, NC 28403 }}{\footnotesize \par}

\medskip{}

\begin{abstract}
We study the traveling wave solutions to a reaction diffusion system
modeling the public goods game with altruistic behaviors. The existence
of the waves is derived through monotone iteration of a pair of classical
upper- and lower solutions. The waves are shown to be unique and strictly
monotonic. A similar KPP wave like asymptotic behaviors are obtained
by comparison principle and exponential dichotomy. The stability of
the traveling waves with non-critical speed is investigated by spectral
analysis in the weighted Banach spaces. 
\end{abstract}
\bigskip{}

\section{\textbf{Introduction \label{Sec1}}}

\setcounter{equation}{0}

Due to their widespread applications in biology, physics and chemistry,
there has been an increasing interest in the traveling wave solutions
to reaction diffusion systems. A fruitful methods have been developed
in deriving the traveling wave solutions, among which the monotone
iteration method is proved to be rather effective. Such method reduces
the existence problem to that of an ordered pair of upper and lower
solutions. However, it is not an easy task to construct the upper
and lower solution pairs in the classical sense, therefore, \cite{04-BoumenirNguyen,19-Ma,25-WangLiRuan,26-WuZou}
turn to search for the generalized upper and lower solution pairs
and have achieved the desired results. Other methods such as geometric
singular perturbation \cite{02-AlexanderGardnerJones,08-Hosono},
phase plane analysis \cite{01-AiChowYi,06-FeiCarr,11-HouLiMeyer,13-Kanel,14-KanelZhou,15-Kan-on,16-Kan-on,18-KolmogorovPetrovskiiPiskunov,21-Sattinger,22-TangFife,23-Volpert,28-XuZhao}
and the fixed point Theorem \cite{19-Ma} are also successfully applied
to various reaction diffusion systems to obtain the wave solutions.

One interesting yet elusive question is the long term behavior of
the traveling wave solutions. In fact, even for the relatively well
studied Lotka Volterra systems less results (see \cite{01-AiChowYi},
\cite{12-Huang}) are known on the asymptotic decay (growth) rates
of the waves and their stability. This question seems to be related
to environmental selection in a given ecosystem. Roughly speaking,
traveling wave is formed by intersecting the unstable subspace from
one equilibrium and stable subspace from another. The stable and unstable
spaces are generally multi-dimensional, therefore, the traveling wave
selects an unstable direction from the unstable equilibrium to emanate,
and then a stable direction from the stable equilibrium to enter.
The selection affects the asymptotic behaviors and stability of the
wave solutions. 

We try to understand those problems through the study of the traveling
wave solutions of the following reaction diffusion system:

\begin{equation}
\left\{ \begin{array}{l}
\frac{\partial\hat{u}}{\partial t}=\frac{\partial^{2}\hat{u}}{\partial x^{2}}+\hat{u}(r_{1}-\frac{\hat{u}+\hat{v}}{k(\hat{u})}-\alpha),\\
\\\frac{\partial\hat{v}}{\partial t}=\frac{\partial^{2}\hat{v}}{\partial x^{2}}+\hat{v}(r_{2}-\frac{\hat{u}+\hat{v}}{k(\hat{u})}),\\
\\u(x,0)=\phi(x),\; v(x,0)=\psi(x)\end{array}\qquad(x,t)\in\mathbb{R}\times\mathbb{R}^{+},\right.\label{eq: 1.01}\end{equation}
where $\hat{u}=\hat{u}(x,t),\;\hat{v}=\hat{u}(x,t)$ are two competing
populations and $r_{1}$, $r_{2}$, $\alpha$ are positive constants.
The system is a continuous spatial-temporal version of a public goods
games (\cite{24-Wakano}) and describes the interaction between the
populations $\hat{u}(x,t)$ and $\hat{v}(x,t)$. In the game, the
population $\hat{u}$ employs an altruistic strategy in order to ensure
the survival of both populations; while the population $\hat{v}$'s
strategy is to maximize its own gain. The function $k(\hat{u})=k_{0}+k\hat{u}$
with $k,\: k_{0}>0$ represents the public goods, contributed by the
population $\hat{u}$ and shared by the population $\hat{v}$. In
this sense, we call population $\hat{u}$ the cooperators and population
$\hat{v}$ the defectors. The cooperators have to pay certain penalty
measured by $\alpha\hat{u}$ for such altruistic strategy. As the
outcome of the game, the population $\hat{v}$ eventually wins while
the population $\hat{u}$ loses. So two equilibrium states are of
interest to us - the one dominated by $\hat{u}$ and the one dominated
by $\hat{v}$. The game is played by finding a passage between the
two equilibrium states.

After rescaling the density functions (see \cite{09-HouFengLu}),
we may further assume that $r_{1}=r_{2}=k_{0}=1$. The rescaled system
has two unstable constant equilibrium states $(0,0)$, $(K^{*},0)$
with $K^{*}=\frac{1-\alpha}{1-k+\alpha k}>0$ and one asymptotically
stable state $(0,1)$. Throughout the rest of the paper we make the
following assumption:

\medskip{}
$\qquad$ \textbf{H}: $0<\alpha,\: k<1.$

\medskip{}

The transformations $\hat{u}=K^{*}-u,\;\hat{v}=v$ change \eqref{eq: 1.01}
into a local monotone type. Disregarding the initial conditions, the
traveling wave solution to \eqref{eq: 1.01} has the form $(u(\xi),v(\xi))^{T}:=(u(x+ct),v(x+ct))^{T}$,
$c>0$ and solves the following system \begin{equation}
\left\{ \begin{array}{l}
u''-cu'-(K^{*}-u)(1-\frac{k^{*}-u+v}{1+k(K^{*}-u)}-\alpha)=0,\\
\\v''-cv'+v(1-\frac{K^{*}-u+v}{1+k(K^{*}-u)})=0,\\
\\(u,v)^{T}(-\infty)=(0,0)^{T},\;(u,v)^{T}(+\infty)=(K^{*},1)^{T}.\end{array}\right.\label{eq:1.02}\end{equation}
where $T$ represents the vector transposition and the prime the derivative
with respect to $\xi$.

In our derivation of the traveling wave solutions, the classical results
on KPP (Kolmogorov-Petrovskii-Piskunov, \cite{18-KolmogorovPetrovskiiPiskunov})
equation and the monotone structure of \eqref{eq:1.02} play a central
role. Noting that if $u$ is replaced by zero then the second equation
in \eqref{eq:1.02} is a KPP equation. The construction of the upper-
and lower solution pairs is based on this simple observation. Cares
also be taken on the asymptotic decay rates of the upper and lower
solution pairs to ensure that they are ordered and have the same decay
rate at negative infinity. The existence of the traveling wave solutions
for system \eqref{eq:1.02} then follows from the monotone iteration
scheme by \cite{04-BoumenirNguyen} and \cite{26-WuZou}, and the
asymptotics of the traveling wave will be obtained from the exponential
dichotomy and comparison. 

The stability of the traveling wave solutions with non-critical speed
is investigated by spectral analysis. Since the wave solutions are
essentially unstable in the unweighted spaces such as $L^{\infty}$
or $C_{\mbox{unif}}$, similar to the stability issues of KPP waves
\cite{10-HouLi,21-Sattinger,27-WuLi,29-WuXingYe}, we first localize
the initial perturbations to certain exponentially weighted Banach
spaces. Then following the standard treatment of such problem as in
\cite{02-AlexanderGardnerJones,17-Kapitula,20-Sandstede,21-Sattinger},
we show that the linearized operator around the wave solution, with
noncritical wave speeds, does not have eigenvalues with positive or
zero real part. Several new ideas are incorporated to locate the eigenvalues
of the linearized operator. This, along with the fact that the linearized
operator is sectorial leads to the nonlinear stability of the wave
solutions. The stability of the wave solution implies that a sufficiently
small initial disturbance will not change the outcome of the game.
We remark that though the paper is written for a particular system
\eqref{eq: 1.01}, the method presents here can be easily adapted
to other reaction diffusion systems with monstable structure. We further
remark that the stability of the waves with critical wave speed has
a complete different nature and requires a different treatment, so
we leave it open for now.

The paper is organized as follows: in Section 2 we show the existence,
uniqueness, monotonicity, and derive the asymptotic decay (growth)
rates of the traveling wave solutions. In Section 3 we study the stability
of the traveling wave solutions in weighted Banach spaces.

\section{\textbf{Existence and Asymptotic Decay Rates \label{Sec2}}}

\setcounter{equation}{0}

In this section, we establish the existence of traveling wave solutions
and derive their asymptotic decay rates through monotone iteration
of a pair of smooth upper and lower solutions. The definition of the
upper and lower solutions is standard. For the rest of the paper the
inequality between two vectors is component-wise.
\begin{defn}
A smooth function $(\bar{u}(\xi),\bar{v}(\xi))^{T}$, $\xi\in\mathbb{R}$
is an upper solution of \eqref{eq:1.02} if its derivatives $(\bar{u}',\bar{v}')^{T}$
and $(\bar{u}'',\bar{v}'')^{T}$ are continuous on $\mathbb{R}$,
and $(\bar{u},\bar{v})^{T}$ satisfies \begin{equation}
\left\{ \begin{array}{l}
u''-cu'-(K^{*}-u)(1-\frac{K^{*}-u+v}{1+k(K^{*}-u)}-\alpha)\leq0,\\
\\v''-cv'+v(1-\frac{K^{*}-u+v}{1+k(K^{*}-u)})\leq0,\end{array}\right.\label{eq:2.01}\end{equation}
with the boundary conditions \begin{equation}
\left(\begin{array}{c}
u\\
v\end{array}\right)(-\infty)=\left(\begin{array}{c}
0\\
0\end{array}\right),\,\;\left(\begin{array}{c}
u\\
v\end{array}\right)(+\infty)\geq\left(\begin{array}{c}
K^{*}\\
1\end{array}\right).\label{eq:2.02}\end{equation}

A lower solution of \eqref{eq:1.02} is defined in a similar way by
reversing the inequalities in \eqref{eq:2.01} and \eqref{eq:2.02}. 
\end{defn}
The construction of the smooth upper and lower solution pairs is based
on the solutions of the KPP equation: \begin{equation}
\left\{ \begin{array}{l}
w''-cw'+f(w)=0,\\
\\w(-\infty)=0,\quad w(+\infty)=b.\end{array}\right.\label{eq:2.03}\end{equation}
where $f\in C^{2}([0,\, b])$ and $f>0$ on the open interval $(0,b)$
with $f(0)=f(b)=0$, $f'(0)=a_{1}>0$ and $f'(b)=-b_{1}<0$. We first
recall the following result (\cite{21-Sattinger}):
\begin{lem}
\label{lem:2.02}Corresponding to every $c\geq2\sqrt{a_{1}}$, system
\eqref{eq:2.03} has a unique (up to a translation of the origin)
monotonically increasing traveling wave solution $w(\xi)$ for $\xi\in\mathbb{R}$.
The traveling wave solution $w$ has the following asymptotic behaviors:

For the wave solution with non-critical speed $c>2\sqrt{a_{1}}$,
we have \begin{equation}
w(\xi)=a_{w}e^{\frac{c-\sqrt{c^{2}-4a_{1}}}{2}\xi}+o(e^{\frac{c-\sqrt{c^{2}-4a_{1}}}{2}\xi})\mbox{ as }\xi\rightarrow-\infty,\label{eq:2.04}\end{equation}
 \begin{equation}
w(\xi)=b-b_{w}e^{\frac{c-\sqrt{c^{2}+4b_{1}}}{2}\xi}+o(e^{\frac{c-\sqrt{c^{2}+4b_{1}}}{2}\xi})\mbox{ as }\xi\rightarrow+\infty,\label{eq:2.05}\end{equation}
 where $a_{w}$ and $b_{w}$ are positive constants.

For the wave with critical speed $c=2\sqrt{a_{1}}$, we have

\[
w(\xi)=(a_{c}+d_{c}\xi)e^{\sqrt{a_{1}}\xi}+o(\xi e^{\sqrt{a_{1}}\xi})\mbox{ as }\xi\rightarrow-\infty,\]

\[
w(\xi)=b-b_{c}e^{(\sqrt{a_{1}}-\sqrt{a_{1}+b_{1}})\xi}+o(e^{(\sqrt{a_{1}}-\sqrt{a_{1}+b_{1}})\xi})\mbox{ as }\xi\rightarrow+\infty,\]
 where the constant $d_{c}$ is negative, $b_{c}$ is positive and
$a_{c}\in\mathbb{R}$.
\end{lem}
To construct the upper-solution for the system \eqref{eq:1.02}, we
begin with the following form of KPP system: \begin{equation}
\left\{ \begin{array}{l}
\tilde{v}''-c\tilde{v}'+{\displaystyle \frac{\alpha}{(1-k+\alpha k)(1+kK^{*}(1-\tilde{v}))}\tilde{v}(1-\tilde{v})}=0,\\
\\\tilde{v}(-\infty)=0,\quad\tilde{v}(+\infty)=1,\end{array}\right.\label{eq:2.06}\end{equation}
where relating to \eqref{eq:2.03}, $f(\tilde{v})=\frac{\alpha}{(1-k+\alpha k)(1+kK^{*}(1-\tilde{v}))}\tilde{v}(1-\tilde{v})>0$
for $\tilde{v}\in(0,1)$. $f(0)=f(1)=0$, $f'(0)=\alpha>0$ and $f'(1)=-\alpha/(1-k+\alpha k)<0$.
According to Lemma \ref{lem:2.02}, for each fixed $c\geq2\sqrt{\alpha}$,
system \eqref{eq:2.06} has a unique (up to a translation of the origin)
traveling wave solution $\bar{v}(\xi)$ satisfying the given boundary
conditions. Define \begin{equation}
\left(\begin{array}{c}
\bar{u}(\xi)\\
\\\bar{v}(\xi)\end{array}\right)=\left(\begin{array}{c}
K^{*}\bar{v}(\xi)\\
\\\bar{v}(\xi)\end{array}\right),\quad\xi\in\mathbb{R},\label{eq:2.07}\end{equation}
 then we have the following result,
\begin{lem}
\label{lem:2.03}For each fixed $c\geq2\sqrt{\alpha}$, \eqref{eq:2.07}
is a smooth upper solution for system \eqref{eq:1.02}. \end{lem}
\begin{proof}
On the boundary, one has $(\bar{u},\bar{v})^{T}(-\infty)=(0,0)^{T}$,
$(\bar{u},\bar{v})^{T}(+\infty)=(K^{*},\,1)^{T}$. 

As for the $u$ component, we have

\[
\begin{array}{ll}
 & \bar{u}''-c\bar{u}'-(K^{*}-\bar{u})(1-\alpha-{\displaystyle \frac{K^{*}-K^{*}\bar{v}+\bar{v}}{1+k(K^{*}-\bar{u})}})\\
\\= & K^{*}[\bar{v}''-c\bar{v}'-(1-\bar{v})(1-\alpha-{\displaystyle \frac{K^{*}-K^{*}\bar{v}+\bar{v}}{1+kK^{*}(1-\bar{v})})}]\\
\\= & {\displaystyle -(1-\bar{v})[\frac{\frac{\alpha}{1-k+\alpha k}\bar{v}-K^{*}+K^{*}\bar{v}-\bar{v}+1+kK^{*}(1-\bar{v})}{1+kK^{*}(1-\bar{v})}-\alpha]}\\
\\= & -(1-\bar{v}){\displaystyle \frac{\alpha kK^{*}\bar{v}}{1+kK^{*}(1-\bar{v})}}\leq0.\end{array}\]

We then verify the second component in \eqref{eq:1.02}, \[
\begin{array}{ll}
 & \bar{v}''-c\bar{v}'+\bar{v}[1-{\displaystyle \frac{K^{*}-K^{*}\bar{v}+\bar{v}}{1+kK^{*}(1-\bar{v})}}]\\
\\= & {\displaystyle -\frac{\alpha}{(1+kK^{*}(1-\tilde{v}))(1-k+\alpha k)}}\bar{v}(1-\bar{v})+\bar{v}[1-{\displaystyle \frac{K^{*}-K^{*}\bar{v}+\bar{v}}{1+kK^{*}(1-\bar{v})}}]\\
\\= & {\displaystyle \bar{v}[\frac{\frac{\alpha}{1-k+\alpha k}-\bar{v}(1+kK^{*}-K^{*})-\frac{\alpha}{1-k+\alpha k}(1-\bar{v})}{1+kK^{*}(1-\bar{v})}]}=0,\end{array}\]
since $1+kK^{*}-K^{*}=\frac{\alpha}{1-k+\alpha k}$.

Thus $(\bar{u},\bar{v})$ forms a smooth upper-solution for (\ref{eq:1.02}). 
\end{proof}
We next construct the lower solution pair for system (\ref{eq:1.02}).
For a small but fixed number $l$ with $0<l<1-k+k\alpha$, we begin
with yet another KPP system: \begin{equation}
\left\{ \begin{array}{l}
\check{v}''-c\check{v}'+{\displaystyle \frac{\alpha}{(1-k+\alpha k)(1+kK^{*}(1-l\check{v}))}\check{v}(1-\frac{1+kK^{*}-lK^{*}}{1+kK^{*}-K^{*}}\check{v})}=0,\\
\\\check{v}(-\infty)=0,\quad\check{v}(+\infty)={\displaystyle \frac{1+kK^{*}-K^{*}}{1+kK^{*}-lK^{*}}}<1.\end{array}\right.\label{eq:2.08}\end{equation}
 Corresponding to the notions in Lemma \ref{lem:2.02}, \[
f(\check{v})={\displaystyle \frac{\alpha}{(1-k+\alpha k)(1+kK^{*}(1-l\check{v}))}}\check{v}(1-\frac{1+kK^{*}-lK^{*}}{1+kK^{*}-K^{*}}\check{v})>0\]
 for $\check{v}\in(0,\frac{1+kK^{*}-K^{*}}{1+kK^{*}-lK^{*}})$. $f(0)=f(\frac{1+kK^{*}-K^{*}}{1+kK^{*}-lK^{*}})=0$,
$f'(0)=\alpha>0$, and $f'(\frac{1+kK^{*}-K^{*}}{1+kK^{*}-lK^{*}})=-\frac{1-l+l\alpha}{1-l+l\alpha(1-k+\alpha k)}<0$. 

For each fixed $c\geq2\sqrt{\alpha}$, define \begin{equation}
\left(\begin{array}{c}
\underline{u}(\xi)\\
\\\underline{v}(\xi)\end{array}\right)=\left(\begin{array}{c}
K^{*}l\underline{v}(\xi)\\
\\\underline{v}(\xi)\end{array}\right),\quad\xi\in\mathbb{R}\label{eq:2.09}\end{equation}
 with $\underline{v}(\xi)$ a solution of (\ref{eq:2.08}).
\begin{lem}
\label{lem:2.04}For each $c\geq2\sqrt{\alpha}$ , \eqref{eq:2.09}
is a smooth lower solution of system (\ref{eq:1.02}). \end{lem}
\begin{proof}
On the boundary, one has \[
\left(\begin{array}{c}
\underline{u}(-\infty)\\
\\\underline{v}(-\infty)\end{array}\right)=\left(\begin{array}{c}
0\\
\\0\end{array}\right),\]
 and\[
\left(\begin{array}{c}
\underline{u}(+\infty)\\
\\\underline{v}(+\infty)\end{array}\right)=\left(\begin{array}{c}
{\displaystyle K^{*}l\frac{1+kK^{*}-K^{*}}{1+kK^{*}-lK^{*}}}\\
\\{\displaystyle \frac{1+kK^{*}-K^{*}}{1+kK^{*}-lK^{*}}}\end{array}\right)<\left(\begin{array}{c}
K^{*}\\
\\1\end{array}\right).\]

Furthermore,\[
\begin{array}{ll}
 & {\displaystyle \underline{u}''-c\underline{u}'-(K^{*}-\underline{u})(1-\alpha-\frac{K^{*}-\underline{u}+\underline{v}}{1+k(K^{*}-\underline{u})})}\\
\\= & {\displaystyle K^{*}l[\underline{v}''-c\underline{v}'-(\frac{1}{l}-\underline{v})(1-\alpha-\frac{K^{*}-\underline{u}+\underline{v}}{1+k(K^{*}-\underline{u})})}]\\
\\= & {\displaystyle K^{*}l[-(\frac{1}{l}-\underline{v})(1-\alpha-\frac{K^{*}-\underline{u}+\underline{v}}{1+k(K^{*}-\underline{u})})}\\
\\ & {\displaystyle -\frac{\alpha}{(1+kK^{*}(1-l\underline{v}))(1-k+\alpha k)}\underline{v}(1-\frac{1+kK^{*}-lK^{*}}{1+kK^{*}-K^{*}}\underline{v})]}\\
\\= & {\displaystyle -K^{*}l\{(\frac{1}{l}-\underline{v})\frac{(1-\alpha)(1+k(K^{*}-\underline{u}))-K^{*}+\underline{u}-\underline{v}}{1+k(K^{*}-\underline{u})}}\\
\\ & {\displaystyle +\frac{\underline{v}[1+kK^{*}-K^{*}-(1+kK^{*}-lK^{*})\underline{v}]}{1+k(K^{*}-\underline{u})}\}}\\
\\= & {\displaystyle -\frac{K^{*}l}{1+k(K^{*}-\underline{u})}\underline{v}\{(\frac{1}{l}-\underline{v})[K^{*}l-1-k(1-\alpha)K^{*}l]+\frac{\alpha}{1-k+\alpha k}-\frac{1-l(1-\alpha)}{1-k+\alpha k}\underline{v}\}}\\
\\= & -{\displaystyle \frac{K^{*2}l}{1+k(K^{*}-\underline{u})}}\underline{v}[1-k+\alpha k-\frac{1-k+\alpha k}{l(1-\alpha)}+\frac{\alpha}{1-\alpha}-k\underline{v}+k(1-\alpha)l\underline{v}]\\
\\\geq & {\displaystyle 0},\end{array}\]
since $0<l<1-k+k\alpha<\frac{1-k(1-\alpha)}{1-k(1-\alpha)^{2}}$ and 

\[
\begin{array}{ll}
 & {\displaystyle \underline{v}''-c\underline{v}'+\underline{v}(1-\frac{K^{*}-\underline{u}+\underline{v}}{1+k(K^{*}-\underline{u})})}\\
\\= & {\displaystyle -\frac{\alpha}{(1+kK^{*}(1-l\underline{v}))(1-k+\alpha k)}\underline{v}(1-\frac{1+kK^{*}-lK^{*}}{1+kK^{*}-K^{*}}\underline{v})+\underline{v}(1-\frac{K^{*}-\underline{u}+\underline{v}}{1+k(K^{*}-\underline{u})})}\\
\\= & {\displaystyle \underline{v}[1-\frac{K^{*}-K^{*}l\underline{v}+\underline{v}}{1+k(K^{*}-\underline{u})}-\frac{\alpha}{(1+kK^{*}(1-\underline{lv}))(1-k+\alpha k)}(1-\frac{1+kK^{*}-lK^{*}}{1+kK^{*}-K^{*}}\underline{v})]}\\
\\= & {\displaystyle \underline{v}[\frac{1+kK^{*}(1-l\underline{v})-K^{*}+K^{*}l\underline{v}-\underline{v}-\frac{\alpha}{1-k+\alpha k}(1-\frac{1+kK^{*}-lK^{*}}{1+kK^{*}-K^{*}}\underline{v})}{1+kK^{*}(1-l\underline{v})}]}\\
\\= & {\displaystyle \underline{v}[1-\frac{1+kK^{*}-K^{*}-\frac{\alpha}{1-k+\alpha k}-\underline{v}(1+klK^{*}+K^{*}l-(1+kK^{*}-lK^{*}))}{1+kK^{*}(1-l\underline{v})})}\\
\\= & {\displaystyle \underline{v}^{2}\frac{kK^{*}(1-l)}{1+kK^{*}(1-l\underline{v})}\geq0}\end{array}\]
 Therefore, the conclusion of the lemma follows. 
\end{proof}
We next show, by shifting the upper solution far to the left, that
the upper- and the lower-solutions as derived in Lemma \ref{lem:2.03}
and Lemma \ref{lem:2.04} are ordered.
\begin{lem}
\label{lem:2.05}Let $c\geq2\sqrt{\alpha}$ be fixed and $(\bar{u},\bar{v})^{T}$,
$(\underline{u},\underline{v})^{T}$ be respectively the upper and
lower solutions defined in \eqref{eq:2.07} and \eqref{eq:2.09},
then there exists a number $r\geq0$, such that $(\bar{u},\bar{v})^{T}(\xi+r)\geq(\underline{u},\underline{v})^{T}(\xi)$
for $\xi\in\mathbb{R}.$ \vspace{5pt}
\end{lem}
\begin{proof}
The proof is only for $c>2\sqrt{\alpha}$ since the one for $c=2\sqrt{\alpha}$
is similar. We first derive the asymptotic behaviors of the upper-
and lower-solutions at infinities. By Lemma \ref{lem:2.02}, \begin{equation}
\left(\begin{array}{c}
\bar{u}\\
\\\bar{v}\end{array}\right)(\xi)=\left(\begin{array}{c}
K^{*}A_{1}\\
\\A_{1}\end{array}\right)e^{\frac{c-\sqrt{c^{2}-4\alpha}}{2}\xi}+o(e^{\frac{c-\sqrt{c^{2}-4\alpha}}{2}\xi})\label{eq:2.10}\end{equation}
 and

\begin{equation}
\left(\begin{array}{c}
\underline{u}\\
\\\underline{v}\end{array}\right)(\xi)=\left(\begin{array}{c}
K^{*}lB_{1}\\
\\B_{1}\end{array}\right)e^{\frac{c-\sqrt{c^{2}-4\alpha}}{2}\xi}+o(e^{\frac{c-\sqrt{c^{2}-4\alpha}}{2}\xi})\label{eq:2.11}\end{equation}
 as $\xi\rightarrow-\infty$; and

\begin{equation}
\left(\begin{array}{c}
\bar{u}\\
\\\bar{v}\end{array}\right)(\xi)=\left(\begin{array}{c}
K^{*}\\
\\1\end{array}\right)-\left(\begin{array}{c}
K^{*}A_{2}\\
\\A_{2}\end{array}\right)e^{\frac{c-\sqrt{c^{2}+4\alpha/(1-k+\alpha k)}}{2}\xi}+o(e^{\frac{c-\sqrt{c^{2}+4\alpha}}{2}\xi})\label{eq:2.12}\end{equation}
 and

\begin{equation}
\left(\begin{array}{c}
\underline{u}\\
\\\underline{v}\end{array}\right)(\xi)=\left(\begin{array}{c}
K^{*}l\\
\\1\end{array}\right)\frac{1+kK^{*}-K^{*}}{1+kK^{*}-lK^{*}}-\left(\begin{array}{c}
K^{*}lB_{2}\\
\\B_{2}\end{array}\right)e^{\underline{\mu}\xi}+o(e^{\underline{\mu}\xi})\label{eq:2.13}\end{equation}
 as $\xi\rightarrow+\infty$, where $\underline{\mu}=\frac{1}{2}(c-\sqrt{c^{2}+4\frac{1-l+l\alpha}{1-l+l\alpha(1-k+\alpha k)}})<0$
and $A_{1}$, $A_{2}$, $B_{1}$, \textbf{$B_{2}$} are positive constants.

Since (\ref{eq:2.06}) is translation invariant, $\bar{v}^{\tilde{r}}(\xi)\equiv\bar{v}(\xi+\tilde{r})$
is also a solution of (\ref{eq:2.06}) for any $\tilde{r}\in\mathbb{R}$.
It then follows that $(\bar{u}^{\tilde{r}},\bar{v}^{\tilde{r}})^{T}(\xi)$
is an upper solution for system (\ref{eq:1.02}). For the asymptotic
behavior of $(\bar{u},\bar{v})^{\tilde{r}}(\xi)$ at $-\infty$, we
can simply replace $(K^{*}A_{1},\: A_{1})^{T}$ by $(K^{*}A_{1},\: A_{1})^{T}e^{\frac{1}{2}(c-\sqrt{c^{2}-4\alpha})\,\tilde{r}}$
in (\ref{eq:2.10}). Now we choose $\tilde{r}>0$ large enough such
that

\[
(K^{*}A_{1},\: A_{1})^{T}\, e^{\frac{c-\sqrt{c^{2}-4\alpha}}{2}\tilde{r}}>(K^{*}A_{1},\: A_{1})^{T}.\]
 Then there exists a sufficiently large $N_{1}>0$ such that

\begin{equation}
\left(\begin{array}{c}
\bar{u}^{\tilde{r}}(\xi)\\
\\\bar{v}^{\tilde{r}}(\xi)\end{array}\right)>\left(\begin{array}{c}
\underline{u}(\xi)\\
\\\underline{v}(\xi)\end{array}\right)\quad{\normalcolor \mbox{for}\;}\xi\in(-\infty,-N_{1}].\label{eq:2.14}\end{equation}

On the other hand, the boundary conditions of the upper and lower
solutions at $+\infty$ also imply that there exists a number $N_{2}>0$
such that

\begin{equation}
\left(\begin{array}{c}
\bar{u}^{\tilde{r}}(\xi)\\
\\\bar{v}^{\tilde{r}}(\xi)\end{array}\right)>\left(\begin{array}{c}
\underline{u}(\xi)\\
\\\underline{v}(\xi)\end{array}\right)\quad{\normalcolor \mbox{for}\;}\xi\in[N_{2},\:+\infty).\label{eq:2.15}\end{equation}

We next show that the inequalities (\ref{eq:2.14}) and (\ref{eq:2.15})
also hold on the interval $[-N_{1},N_{2}]$. There are two possible
cases to deal with:

Case 1. If we already have

\begin{equation}
\left(\begin{array}{c}
\bar{u}^{\tilde{r}}(\xi)\\
\\\bar{v}^{\tilde{r}}(\xi)\end{array}\right)\geq\left(\begin{array}{c}
\underline{u}(\xi)\\
\\\underline{v}(\xi)\end{array}\right)\quad\mbox{on }[-N_{1},N_{2}],\label{eq:2.16}\end{equation}
 we then let $\tilde{r}=r$ and obtain the Lemma.

Case 2. There exists a point $\xi_{0}\in(-N_{1},\: N_{2})$ such that

\begin{equation}
\left(\begin{array}{c}
\bar{u}^{\tilde{r}}(\xi_{0})\\
\\\bar{v}^{\tilde{r}}(\xi_{0})\end{array}\right)\leq\left(\begin{array}{c}
\underline{u}(\xi_{0})\\
\\\underline{v}(\xi_{0})\end{array}\right).\label{eq:2.17}\end{equation}
 with strict inequality holding for at least one of the two components.

In this case, we will use the the \textit{Sliding Domain method}.
We first shift $(\bar{u}^{\tilde{r}},\bar{v}^{\tilde{r}})$ to the
left by increasing $\tilde{r}$ until we can find a $r_{1}>\tilde{r}>0$
such that $(\bar{u}^{r_{1}}(\xi),\bar{v}^{r_{1}}(\xi))^{T}>(\underline{u}(\xi),\underline{v}(\xi))^{T}$
on the interval $[-N_{1},\, N_{2}-(r_{1}-\tilde{r})]$. We then shift
$(\bar{u}^{r_{1}}(\xi),\bar{v}^{r_{1}}(\xi))^{T}$ back to the right
by decreasing $r_{1}$ to some $r_{2}>\tilde{r}$ such that one of
the branches of the upper solution touches its counterpart of the
lower solution at some point $\xi_{2}$ in the interval $(-N_{1}+r_{2},\, N_{2}-(r_{1}-\tilde{r}))$.
On the endpoints of the interval $(-N_{1}+r_{2},N_{2}-(r_{1}-\tilde{r}))$,
we still have $(\bar{u}^{r_{2}}(\xi),\bar{v}^{r_{2}}(\xi))^{T}>(\underline{u}(\xi),\underline{v}(\xi))^{T}$.
In summary, we now have $\bar{u}^{r_{2}}(\xi_{2})=\underline{u}(\xi_{2})$
and $\bar{u}^{r_{2}}(\xi)\geq\underline{u}(\xi)$, $\bar{v}^{r_{2}}(\xi)\geq\underline{v}(\xi)$
for $\xi\in(-N_{1}+r_{2},N_{2}-(r_{1}-\tilde{r}))$.

Letting $W(\xi):=(\bar{u}^{r_{2}},\bar{v}^{r_{2}})^{T}(\xi)-(\underline{u},\underline{v})^{T}(\xi)$
and $F=(F_{1},F_{2})^{T}=(-(K^{*}-u)(1-\alpha-\frac{K^{*}-u+v}{1+k(K^{*}-u)}),v(1-\frac{K^{*}-u+v}{1+k(K^{*}-u)}))^{T}$.
For $\xi\in(-N_{1}+r_{2},N_{2}-(r_{1}-\tilde{r}))$ we have, \begin{equation}
W''-cW'+\left[\begin{array}{cc}
\frac{\partial F_{1}}{\partial u}(\underline{u}+\zeta_{1}w_{1},\bar{v}^{r_{2}}), & \frac{\partial F_{1}}{\partial v}(\bar{u}^{r_{2}},\underline{v}+\zeta_{2}w_{2})\\
\\\frac{\partial F_{2}}{\partial u}(\underline{u}+\zeta_{3}w_{1},\bar{v}^{r_{2}}), & \frac{\partial F_{2}}{\partial v}(\bar{u}^{r_{2}},\underline{v}+\zeta_{4}w_{2})\end{array}\right]W=0,\label{eq:2.18}\end{equation}
 for some $\zeta_{i}\in[0,1]$, $i=1,2,3,4$. Since the above system
is monotone and the cube $[(0,0),(K^{*},1)]$ is convex, we can readily
deduce by Maximum Principle that $W(\xi)>0$ for $\xi\in[-N_{1}+r_{2},N_{2}-(r_{1}-\tilde{r})]$.
This is a contradiction, i.e., such $\xi_{2}$ does not exist and
we can further decrease $r_{2}$ to $\tilde{r}$. This shows that
$\xi_{0}$ does not exist either. We therefore have $(\bar{u}^{\tilde{r}},\bar{v}^{\tilde{r}})^{T}(\xi)\geq(\underline{u},\underline{v})^{T}(\xi)$
for $\xi\in\mathbb{R}$. 
\end{proof}
To ease the burden of notations, we still use $(\bar{u},\bar{v})^{T}$
to denote the shifted upper solution as given in lemma \ref{lem:2.05}.
With such constructed ordered upper and lower solution pairs, we now
have,
\begin{thm}
\label{thm:2.06}For every $c\geq2\sqrt{\alpha}$, system (\ref{eq:1.02})
has correspondingly a unique (up to a translation of the origin) traveling
wave solution. The traveling wave solution is strictly increasing
and has the following asymptotic properties:

1. Corresponding to the wave speed $c>2\sqrt{\alpha}$, \begin{equation}
\left(\begin{array}{c}
u(\xi)\\
\\v(\xi)\end{array}\right)=\left(\begin{array}{c}
A_{1}\\
\\A_{2}\end{array}\right)e^{\frac{c-\sqrt{c^{2}-4\alpha}}{2}\xi}+o(e^{\frac{c-\sqrt{c^{2}-4\alpha}}{2}\xi})\label{eq:2.19}\end{equation}
 as $\xi\rightarrow-\infty$; and \begin{equation}
\left(\begin{array}{c}
u(\xi)\\
\\v(\xi)\end{array}\right)=\left(\begin{array}{c}
K^{*}\\
\\1\end{array}\right)-\left(\begin{array}{c}
\bar{A}_{1}\\
\\\bar{A}_{2}\end{array}\right)e^{\frac{c-\sqrt{c^{2}+4\alpha}}{2}\xi}+o(e^{\frac{c-\sqrt{c^{2}+4\alpha}}{2}\xi})\label{eq:2.20}\end{equation}
 as $\xi\rightarrow+\infty,$ where $A_{1}$, $A_{2}$, $\bar{A}_{1}$,
$\bar{A}_{2}$ are positive constants.

2. Corresponding to the wave speed $c_{\mbox{critical}}=2\sqrt{\alpha}$,
we have

\begin{equation}
\left(\begin{array}{c}
u(\xi)\\
\\v(\xi)\end{array}\right)=\left(\begin{array}{c}
(A_{11c}+A_{12c}\xi)\\
\\(A_{21c}+A_{22c}\xi)\end{array}\right)e^{\sqrt{\alpha}\xi}+o(\xi e^{\sqrt{\alpha}\xi})\label{eq:2.21}\end{equation}
 as $\xi\rightarrow-\infty$; and

\begin{equation}
\left(\begin{array}{c}
u(\xi)\\
\\v(\xi)\end{array}\right)=\left(\begin{array}{c}
K^{*}\\
\\1\end{array}\right)-\left(\begin{array}{c}
\bar{A}_{11}\\
\\\bar{A}_{22}\end{array}\right)e^{(1-\sqrt{2})\sqrt{\alpha}\xi}+o(e^{(1-\sqrt{2})\sqrt{\alpha}\xi})\label{eq:2.22}\end{equation}
 as $\xi\rightarrow+\infty,$ where $A_{12c},\, A_{22c}<0$, $A_{11c}$,
$A_{21c}\in\mathbb{R}$ and $\bar{A}_{11},$$\bar{A}_{22}>0$. \end{thm}
\begin{proof}
Starting from the upper and lower solution pairs obtained in Lemma
\ref{lem:2.03}-\ref{lem:2.04} and using the monotone iteration scheme
given in \cite{04-BoumenirNguyen,26-WuZou}, we obtain the existence
of the solution $(u(\xi),v(\xi))^{T}$ to the first two equations
in \eqref{eq:1.02} for every $c\geq2\sqrt{\alpha}$, which satisfies
\[
\left(\begin{array}{c}
\underline{u}(\xi)\\
\\\underline{v}(\xi)\end{array}\right)\leq\left(\begin{array}{c}
u(\xi)\\
\\v(\xi)\end{array}\right)\leq\left(\begin{array}{c}
\bar{u}(\xi)\\
\\\bar{v}(\xi)\end{array}\right).\]
On the boundary, it is easy to see that the solution tends to $(0,0)^{T}$
as $\xi\rightarrow-\infty$, and $(K^{*},1)^{T}$ as $\xi\rightarrow+\infty$
according to the last inequality.

To derive the asymptotic decay rate of the traveling wave solutions
at $\pm\infty$, we let $c\geq2\sqrt{\alpha}$ and \begin{equation}
U(\xi):=(u(\xi),v(\xi))^{T}\quad\mbox{for}\:-\infty<\xi<\infty\label{eq:2.23}\end{equation}
 be the corresponding traveling wave solution of \eqref{eq:1.02}
generated from the monotone iteration. We differentiate \eqref{eq:1.02}
with respect to $\xi$, and note that $U'(\xi):=(w_{1},w_{2})^{T}(\xi)$
satisfies \begin{equation}
w_{1}''-cw_{1}'+A_{11}(u,v)w_{1}+A_{12}(u,v)w_{2}=0,\label{eq:2.24}\end{equation}
 \begin{equation}
w_{2}''-cw_{2}'+A_{21}(u,v)w_{1}+A_{22}(u,v)w_{2}=0,\label{eq:2.25}\end{equation}
 where \[
\begin{array}{lll}
A_{11}(u,v) & = & 1-{\displaystyle \frac{K^{*}-u+v}{1+k(K^{*}-u)}-\alpha}\\
\\ & + & (K^{*}-u){\displaystyle \frac{kv-1}{(1+k(K^{*}-u))^{2}}},\\
\\A_{12}(u,v) & = & {\displaystyle \frac{K^{*}-u}{1+k(K^{*}-u)}},\\
\\A_{21}(u,v) & = & -{\displaystyle \frac{kv-1}{(1+k(K^{*}-u))^{2}}}v,\\
\\A_{22}(u,v) & = & 1-{\displaystyle \frac{K^{*}-u+v}{1+k(K^{*}-u)}}\\
\\ & - & {\displaystyle \frac{1}{1+k(K^{*}-u)}}v.\end{array}\]

Lemma \ref{lem:2.02} implies that the upper- and the lower-solutions
as derived in Lemma \ref{lem:2.03} and Lemma \ref{lem:2.04} have
the same asymptotic rates at $-\infty$. (\ref{eq:2.19}) and (\ref{eq:2.21})
then follow from Lemma \ref{lem:2.05}.

We next study the exponential decay rates of the traveling wave solution
$U(\xi)$ at $+\infty$. The asymptotic system of (\ref{eq:2.24})
and (\ref{eq:2.25}) as $\xi\rightarrow+\infty$ is \begin{equation}
\left\{ \begin{array}{l}
\psi_{1}''-c\psi_{1}'-\alpha\psi_{1}=0,\\
\\\psi_{2}''-c\psi_{2}'+(1-k)\psi_{1}-\psi_{2}=0.\end{array}\right.\label{eq:2.26}\end{equation}

It is easy to see that system \eqref{eq:2.26} admits exponential
dichotomy \cite{05-Coddington}. Since the traveling wave solution
$(u(\xi),v(\xi))^{T}$ converge monotonically to a constant limit
as $\xi\rightarrow\pm\infty$, the derivative of the traveling wave
solution satisfies $(w_{1}(\pm\infty),w_{2}(\pm\infty))=(0,0)$ (\cite{26-WuZou},
p$658$ Lemma $3.2$). Hence we are only interested in finding bounded
solutions of \eqref{eq:2.26} at $+\infty$.

Since the first equation in \eqref{eq:2.26} is decoupled, we can
write its general solution as \[
\psi_{1}=A^{1}e^{\frac{c+\sqrt{c^{2}+4\alpha}}{2}\xi}+B^{1}e^{\frac{c-\sqrt{c^{2}+4\alpha}}{2}\xi}\]
 for some constants $A^{1}$ and $B^{1}$. Since $w_{1}\rightarrow0$
as $\xi\rightarrow+\infty$, one immediately has $A^{1}=0$.

We then study the second equation of \eqref{eq:2.26}, rewriting the
equation as \begin{equation}
(\psi_{2})''-c(\psi_{2})'-\psi_{2}=-(1-k)\psi_{1},\label{eq:2.27}\end{equation}
 since $\frac{c-\sqrt{c^{2}+4\alpha}}{2}$ is not the characteristics
of the homogeneous equation, one has the following expression for
the general solution of (\ref{eq:2.27}), \begin{equation}
\psi_{2}=\bar{B}_{1}e^{\frac{c-\sqrt{c^{2}+4\alpha}}{2}\xi}+\bar{B}_{2}e^{\frac{c-\sqrt{c^{2}+4}}{2}\xi}+\bar{B}_{3}e^{\frac{c+\sqrt{c^{2}+4}}{2}\xi}.\label{eq:2.28}\end{equation}

Since $w_{2}(\xi)\rightarrow0$ as $\xi\rightarrow+\infty$, then
$\bar{B}_{3}=0.$ Also noticing that (\ref{eq:2.27}) is non-homogeneous,
we have $\bar{B}_{1}\neq0$. By roughness of the exponential dichotomy
(\cite{05-Coddington}) and comparing with the upper solution, we
obtain the asymptotic decay rate of the traveling wave solutions at
$+\infty$ given in (\ref{eq:2.20}) and (\ref{eq:2.23}).

We next show the strict monotonicity of the traveling wave solutions,
which will be a key ingredient in locating the eigenvalues of the
linearized operator in the next section. By the monotone iteration
process\cite{26-WuZou}, the traveling wave solution $U(\xi)$ is
increasing for $\xi\in\mathbb{R}$, it then follows that $(w_{1}(\xi),w_{2}(\xi))^{T}=U'(\xi)\geq0$
and satisfies \begin{equation}
w_{1}''-cw_{1}'+\frac{\partial F_{1}}{\partial u}(u,v)w_{1}+\frac{\partial F_{1}}{\partial v}(u,v)w_{2}=0,\label{eq:2.29}\end{equation}
 \begin{equation}
w_{2}''-cw_{2}'+\frac{\partial F_{2}}{\partial u}(u,v)w_{1}+\frac{\partial F_{2}}{\partial v}(u,v)w_{2}=0,\label{eq:2.30}\end{equation}
 and

\begin{equation}
(w_{1}(\xi),w_{2}(\xi))^{T}\geq0,\;(w_{1},w_{2})^{T}(\pm\infty)=0.\label{eq:2.31}\end{equation}
 The strong Maximum Principle implies that $(w_{1},w_{2})^{T}(\xi)>0$
for $\xi\in\mathbb{R}$. This concludes the strict monotonicity of
the traveling wave solutions.

On the uniqueness of the traveling wave solution for every $c\geq2\sqrt{\alpha}$,
we only prove the conclusion for traveling wave solutions with asymptotics
\eqref{eq:2.19} and \eqref{eq:2.20}, since other case can be proved
similarly. Let $U_{1}(\xi)=(u_{1},v_{1})^{T}$ and $U_{2}(\xi)=(u_{2},v_{2})^{T}$
be two traveling wave solutions of system \eqref{eq: 1.01} with the
same speed $c>2\sqrt{\alpha}$. There exist positive constants $A_{i}$,
$B_{i}$, $i=1,2,3,4$ and a large number $N>0$ such that for $\xi<-N$,\begin{equation}
U_{1}(\xi)=\left(\begin{array}{c}
A_{1}e^{\frac{c+\sqrt{c^{2}-4\alpha}}{2}\xi}\\
\\A_{2}e^{\frac{c+\sqrt{c^{2}-4\alpha}}{2}\xi}\end{array}\right)+o(e^{\frac{c+\sqrt{c^{2}-4\alpha}}{2}\xi})\label{eq:2.32}\end{equation}
 \begin{equation}
U_{2}(\xi)=\left(\begin{array}{c}
A_{3}e^{\frac{c+\sqrt{c^{2}-4\alpha}}{2}\xi}\\
\\A_{4}e^{\frac{c+\sqrt{c^{2}-4\alpha}}{2}\xi}\end{array}\right)+o(e^{\frac{c+\sqrt{c^{2}-4\alpha}}{2}\xi});\label{eq:2.33}\end{equation}
 and for $\xi>N$,\begin{equation}
U_{1}(\xi)=\left(\begin{array}{c}
{\displaystyle K^{*}-B_{1}e^{\frac{c-\sqrt{c^{2}+4\alpha}}{2}\xi}}\\
\\{\displaystyle 1-B_{2}e^{\frac{c-\sqrt{c^{2}+4\alpha}}{2}\xi}}\end{array}\right)+o(e^{\frac{c-\sqrt{c^{2}+4\alpha}}{2}\xi}),\label{eq:2.34}\end{equation}
 \begin{equation}
U_{2}(\xi)=\left(\begin{array}{c}
{\displaystyle K^{*}-B_{3}e^{\frac{c-\sqrt{c^{2}+4\alpha}}{2}\xi}}\\
\\{\displaystyle 1-B_{4}e^{\frac{c-\sqrt{c^{2}+4\alpha}}{2}\xi}}\end{array}\right)+o(e^{\frac{c-\sqrt{c^{2}+4\alpha}}{2}\xi})\label{eq:2.35}\end{equation}
 The traveling wave solutions of system \eqref{eq:1.02} are translation
invariant, thus for any $\theta>0$, $U_{1}^{\theta}(\xi):=U_{1}(\xi+\theta)$
is also a traveling wave solution of \eqref{eq:1.02}. By \eqref{eq:2.32}
and \eqref{eq:2.34}, the solution $U_{1}(\xi+\theta)$ has the asymptotics
\begin{equation}
U_{1}^{\theta}(\xi)=\left(\begin{array}{c}
A_{1}e^{\frac{c+\sqrt{c^{2}-4\alpha}}{2}\theta}e^{\frac{c+\sqrt{c^{2}-4\alpha}}{2}\xi}\\
\\A_{2}e^{\frac{c+\sqrt{c^{2}-4\alpha}}{2}\theta}e^{\frac{c+\sqrt{c^{2}-4\alpha}}{2}\xi}\end{array}\right)+o(e^{\frac{c+\sqrt{c^{2}-4\alpha}}{2}\xi})\label{eq:2.36}\end{equation}
 for $\xi\leq-N$;\begin{equation}
U_{1}^{\theta}(\xi)=\left(\begin{array}{c}
{\displaystyle K^{*}-B_{1}e^{\frac{c-\sqrt{c^{2}+4\alpha}}{2}\theta}e^{\frac{c-\sqrt{c^{2}+4\alpha}}{2}\xi}}\\
\\{\displaystyle 1-B_{2}e^{\frac{c-\sqrt{c^{2}+4\alpha}}{2}\theta}e^{\frac{c-\sqrt{c^{2}+4\alpha}}{2}\xi}}\end{array}\right)+o(e^{\frac{c-\sqrt{c^{2}+4\alpha}}{2}\xi})\label{eq:2.37}\end{equation}
 for $\xi\geq N$.

Choosing $\theta>0$ large enough such that \begin{equation}
A_{1}e^{\frac{c+\sqrt{c^{2}-\alpha}}{2}\theta}>A_{3},\label{eq:2.38}\end{equation}
 \begin{equation}
A_{2}e^{\frac{c+\sqrt{c^{2}-4\alpha}}{2}\theta}>A_{4},\label{eq:2.39}\end{equation}
 \begin{equation}
B_{1}e^{\frac{c-\sqrt{c^{2}+4\alpha}}{2}\theta}<B_{3},\label{eq:2.40}\end{equation}
 \begin{equation}
B_{2}e^{\frac{c-\sqrt{c^{2}+4\alpha}}{2}\theta}<B_{4}.\label{eq:2.41}\end{equation}
then one has \begin{equation}
U_{1}^{\theta}(\xi)>U_{2}(\xi)\label{eq:2.42}\end{equation}
for $\xi\in(-\infty,-N]$$\cup$$[N+\infty).$ We now consider system
\eqref{eq:1.02} on $[-N,+N]$. First, suppose $U_{1}^{\theta}(\xi)\geq U_{2}(\xi)$
on $[-N,+N]$, then the function $W(\xi)=(w_{1}(\xi),w_{2}(\xi))^{T}:=U_{1}^{\theta}(\xi)-U_{2}(\xi)\geq0$
and satisfies for some $\zeta_{i}\in(0,1)$, $i=1,2,3,4$,  \begin{equation}
\begin{array}{ll}
W''-cW'+\left[\begin{array}{cc}
\frac{\partial F_{1}}{\partial u}(u_{2}+\zeta_{1}w_{1},v_{1}), & \frac{\partial F_{1}}{\partial v}(u_{1},v_{2}+\zeta_{2}w_{2})\\
\\\frac{\partial F_{2}}{\partial u}(u_{2}+\zeta_{3}w_{1},v_{1}), & \frac{\partial F_{2}}{\partial v}(u_{1},v_{2}+\zeta_{4}w_{2})\end{array}\right]W=0,\, & \xi\in(-N,N),\\
\\W(-N)>0,\,\,\, W(+N)>0.\end{array}\label{eq:2.43}\end{equation}
Since the above system is monotone, we can readily deduce by the Maximum
Principle that $W(\xi)>0$ on $[-N,N]$. Consequently, we have $U_{1}^{\theta}(\xi)>U_{2}(\xi)$
on $\mathbb{R}$ in this case.

Second, suppose there is some point in $(-N,N)$ such that one of
the components, say the $j$-th component, satisfies $(U_{1}^{\theta}(\xi))_{j}<(U_{2}(\xi))_{j}$
at that point, $j=1$ or $2$. We then increase $\theta$, that is
shift $U_{1}^{\theta}(\xi)$ further left, so that $U_{1}^{\theta}(-N)>U_{2}(-N)$,
$U_{1}^{\theta}(N)>U_{2}(N)$. By the monotonicity of $U_{1}^{\theta}$
and $U_{2}$, we can find a $\bar{\theta}\in(0,2N)$ such that in
the interval $(-N,N)$, we have $U_{1}^{\theta}(\xi+\bar{\theta})>U_{2}(\xi)$.
Shifting $U_{1}^{\theta}(\xi+\bar{\theta})$ back until one component
of $U_{1}^{\theta}(\xi+\bar{\theta})$ first touches its counterpart
of $U_{2}(\xi)$ at some point $\bar{\bar{\xi}}\in[-N,N]$. Since
$U_{1}^{\theta}(\xi+\bar{\theta})$ and $U_{2}(\xi)$ are strictly
increasing, $\bar{\bar{\xi}}\in(-N,N)$, while at $\xi=\pm N$, we
still have $U_{1}^{\theta}(\xi+\bar{\theta})>U_{2}(\xi)$. However,
by Maximum Principle for that component again, we find that component
of $U_{1}^{\theta}$ and $U_{2}$ are identically equal for all $\xi\in[-N,N]$
for a larger $\theta$ than the original one for which \eqref{eq:2.42}
holds. This is a contradiction. Therefore, we must have\[
U_{1}^{\theta}(\xi)>U_{2}(\xi)\]
 for all $\xi\in\mathbb{R}$, where $\theta$ is the one chosen by
means of \eqref{eq:2.38}-\eqref{eq:2.41} as described above.

Now, decrease $\theta$ until one of the following situations happens.

1. There exists a $\bar{\theta}\geq0$, such that $U_{1}^{\bar{\theta}}(\xi)\equiv U_{2}(\xi)$.
In this case we have finished the proof.

2. There exists a $\bar{\theta}\geq0$ and $\xi_{1}\in\mathbb{R}$,
such that one of the components of $U^{\bar{\theta}}$ and $U_{2}$
are equal there; and for all $\xi\in\mathbb{R}$, we have $U_{1}^{\bar{\theta}}(\xi)\geq U_{2}(\xi)$.
On applying the Maximum Principle on $\mathbb{R}$ for that component,
we find $U_{1}^{\bar{\theta}}$ and $U_{2}$ must be identical on
that component. To fix ideas, we suppose that the component is the
first component. Then $U_{1}^{\bar{\theta}}-U_{2}$ satisfies \eqref{eq:2.29}
and \eqref{eq:2.31}. Plugging $w_{1}\equiv0$ into \eqref{eq:2.29}
again we find that there is at least one $\xi_{\bar{\theta}}$ such
that $W_{2}(\xi_{\bar{\theta}})=0$. Then by applying maximum principle
to \eqref{eq:2.30}, we have $w_{2}(\xi)\equiv0$ for $\xi\in\mathbb{R}$.
We have then returned to case 1.

Consequently, in either situation, there exists a $\bar{\theta}\geq0$,
such that \[
U_{1}^{\bar{\theta}}(\xi)\equiv U_{2}(\xi).\]
 for all $\xi\in\mathbb{R}$.
\end{proof}
The next Theorem shows that the lower bound $2\sqrt{\alpha}$ for
the wave speed $c$ is optimal, hence it is the critical minimal wave
speed.
\begin{thm}
\label{thm:2.07}There is no monotone traveling wave solution of \eqref{eq:1.02}
for any $0<c<2\sqrt{\alpha}$. \end{thm}
\begin{proof}
Suppose there is a constant $c$ with $0<c<2\sqrt{\alpha}$ and a
solution $V(\xi)=(v_{1},v_{2})^{T}(\xi)$ of (\ref{eq:1.02}) corresponding
to it. Similar to the proof of Theorem \ref{thm:2.06}, the asymptotic
behaviors of $V(\xi)$ at $-\infty$ are described by \[
\left(\begin{array}{c}
v_{1}(\xi)\\
\\v_{2}(\xi)\end{array}\right)=\left(\begin{array}{c}
A_{s}\\
\\B_{s}\end{array}\right)e^{\frac{c-\sqrt{c^{2}-4\alpha}}{2}\xi}+\left(\begin{array}{c}
\bar{A_{s}}\\
\\\bar{B_{s}}\end{array}\right)e^{\frac{c+\sqrt{c^{2}-4\alpha}}{2}\xi}+h.o.t,\]
 where $(A_{s},B_{s})^{T}$and $(\bar{A_{s},}\bar{B_{s}})$ can not
be both zero due to the non-homogeneity of the limit system of \eqref{eq:2.29}
and \eqref{eq:2.30} at $-\infty$. The condition $0<c<2\sqrt{\alpha}$
implies that $V(\xi)$ is oscillating. This concludes that any solution
of (\ref{eq:1.02}) with $c<2\sqrt{\alpha}$ is not strictly monotone.
 \end{proof}
\begin{rem}
An interesting implication of Theorem \ref{thm:2.06} is that the
populations $u$ and $v$ have the same exponential rate determined
by $\alpha$, the measurement of the penalty for altruism that $u$
receives. The bigger the penalty, the faster the population $u$ decreases
and population $v$ increases. Given the same population growth rate
of $u$ and $v$ , according to the model, the cooperators always
loses the game.
\end{rem}

\section{\textbf{\label{Sec3}Stability of the Traveling Waves}}

\setcounter{equation}{0}

In this section, we investigate the stability of the traveling waves
with non-critical speed. We first show that the traveling wave solution
obtained in Theorem \ref{thm:2.06} with non-critical speed is unstable
in the continuous function space $C(\mathbb{R})\times C(\mathbb{R})$.
This motivates us to investigate the spectrum of the linearized operator
(see \eqref{eq:3.03}) in the exponentially weighted Banach spaces.

Retaining all the transformations and rescalings in section \ref{Sec2},
we consider system \eqref{eq:1.02} in the following equivalent form,
with initial conditions as given in \eqref{eq: 1.01} \[
\left\{ \begin{array}{ccl}
U_{t}(x,t) & = & U_{xx}+F(U),\\
U(x,0) & = & \bar{U(}x).\end{array}\right.\]
 where $U=(u,v)^{T}$ and  $F(U)=(-(K^{*}-u)(1-\alpha-\frac{K^{*}-u+v}{1+k(K^{*}-u)}),v(1-\frac{K^{*}-u+v}{1+k(K^{*}-u)}))^{T}$,
$\bar{U(}x)$ is the initial function.

Rewriting the above equivalent system in terms of $(\xi,t)$ variable,
with moving coordinates $\xi=x+ct$, we have \begin{equation}
\left\{ \begin{array}{ccl}
U_{t}(\xi,t) & = & U_{\xi\xi}-cU_{\xi}+F(U),\\
U(\xi,0) & = & \bar{U}(\xi).\end{array}\right.\label{eq:3.01}\end{equation}

For a fixed wave speed $c^{*}>2\sqrt{\alpha}$ and $\xi=x+c^{*}t$
, $U^{*}(\xi)=(u^{*}(\xi),v^{*}(\xi))^{T}$ obtained in Theorem \ref{thm:2.06}
is then a steady state of \eqref{eq:3.01}. Let $U(\xi,t)=U^{*}(\xi)+V(\xi,t)$
be a solution of \eqref{eq:3.01}. We then obtain the following system
for the perturbation function $V(\xi,t)$:

\begin{equation}
\left\{ \begin{array}{ccl}
V_{t} & = & \mathcal{L}V+\mathcal{R}(V,U^{*}),\\
V(\xi,0) & = & \bar{U}(\xi)-U^{*}(\xi)\end{array}\right.\label{eq:3.02}\end{equation}
 where \begin{equation}
\mathcal{L}V=V_{\xi\xi}-c^{*}V_{\xi}+\frac{\partial F}{\partial U}(U^{*})V\label{eq:3.03}\end{equation}
 is a linear operator, and \begin{equation}
\mathcal{R}(V,U^{*})=F(U^{*}+V)-F(U^{*})-\frac{\partial F}{\partial U}(U^{*})V\label{eq:3.04}\end{equation}
 is a nonlinear operator.

The stability of the traveling wave solution $U^{*}(\xi)$ in a specific
Banach space depends critically on the location of the spectrum $\sigma(\mathcal{L})$
of $\mathcal{L}$. Since $\mathcal{L}$ is defined on $\mathbb{R}$,
it has point spectrum (eigenvalue) \[
\sigma_{p}(\mathcal{L})=\{\lambda\in\sigma(\mathcal{L})\,|\,\lambda\mbox{ is an isolated eigenvalue of finite multiplicity}\}\]
as well as essential spectrum $\sigma_{e}(\mathcal{L})=\sigma(\mathcal{L})\backslash\sigma_{p}(\mathcal{L})$
(\cite{02-AlexanderGardnerJones,17-Kapitula,23-Volpert,21-Sattinger}).
Let $C(\mathbb{R})$ be the space of all continuous functions on $\mathbb{R}$
and $C_{0}$ be a subspace of $C(\mathbb{R})$, \[
C_{0}=\{U(\xi)\,|\: U(\xi)\in C(\mathbb{R})\times C(\mathbb{R}),\,\lim_{|\xi|\rightarrow\infty}||U(\xi)||=0\}\]

\noindent where \[
\parallel U\parallel_{C_{0}}=\sup_{\xi\in\mathbb{R}}\parallel U(\xi)\parallel.\]

\begin{thm}
\textbf{\label{thm:9}} The traveling wave solution $U^{*}(\xi)$
of system (\ref{eq:3.01}) with non-critical wave speed $c^{*}>2\sqrt{\alpha}$
, obtained in Theorem \ref{thm:2.06}, is unstable with initial functions
in $C_{0}$. \vspace{5pt}
\end{thm}
\begin{proof}
We prove the trivial solution of (\ref{eq:3.02}) is unstable. Thus,
it suffices to show that in the space $C_{0\mbox{ }}$ the operator
$\mathcal{L}$ in (\ref{eq:3.03}) has essential spectrum with positive
real parts. It turns out that (\cite{07-Henry,23-Volpert}) the essential
spectrum of the operator $\mathcal{L}$ is bounded by the spectrum
of $\mathcal{L}$ at $\pm\infty$, which are denoted by $\mathcal{L}^{+}$
and $\mathcal{L}^{-}$ respectively. More precisely, we have \begin{equation}
\begin{array}{lll}
\mathcal{L}^{+}V & = & V_{\xi\xi}-c^{*}V_{\xi}+{\displaystyle \frac{\partial F}{\partial U}(U_{+}^{*})V}\\
\\ & = & V_{\xi\xi}-c^{*}V_{\xi}+\left[\begin{array}{cc}
-\alpha & 0\\
\\1-k & -1\end{array}\right]V,\end{array}\label{eq:3.05}\end{equation}
 \begin{equation}
\begin{array}{lll}
\mathcal{L}^{-}V & = & V_{\xi\xi}-c^{*}V_{\xi}+{\displaystyle \frac{\partial F}{\partial U}(U_{-}^{*})V}\\
\\ & = & V_{\xi\xi}-c^{*}V_{\xi}+\left[\begin{array}{cc}
-(1-\alpha)(1-k+\alpha k) & 1-\alpha\\
\\{\displaystyle 0} & \alpha\end{array}\right]V,\end{array}\label{eq:3.06}\end{equation}
 where $U_{\pm}^{*}$ denotes the limit of $U^{*}(\xi)$ as $\xi\rightarrow\pm\infty$.

Now consider the equation \[
\frac{\partial V}{\partial t}=\mathcal{L}^{+}V.\]
 Following \cite{23-Volpert}, we replace $V$ by $e^{(\lambda\tau+i\zeta\xi)}I,$
where $I$ is the identity matrix. We then have \begin{equation}
e^{(\lambda t+i\zeta\xi)}(-\zeta^{2}I-c^{*}\zeta iI+\frac{\partial F}{\partial U}(U_{+}^{*})-\lambda I)=0.\label{eq:3.07}\end{equation}
 The spectrum of the operator $\mathcal{L}^{-}$ consists of curves
\begin{equation}
\det(-\zeta^{2}I-c^{*}\zeta iI+\frac{\partial F}{\partial U}(U_{+}^{*})-\lambda I)=0.\label{eq:3.08}\end{equation}
 Therefore we have \begin{equation}
-\zeta^{2}-c^{*}\zeta i-\alpha-\lambda=0,\label{eq:3.09}\end{equation}
 or \begin{equation}
-\zeta^{2}-c^{*}\zeta i-1-\lambda=0.\label{eq:3.10}\end{equation}

Let $\lambda=x+yi$, then by \eqref{eq:3.09} we have \begin{equation}
x=-\frac{y^{2}}{(c^{*})^{2}}-\alpha,\label{eq:3.11}\end{equation}
 or by \eqref{eq:3.10}, \begin{equation}
x=-\frac{y^{2}}{(c^{*})^{2}}-1.\label{eq:3.12}\end{equation}

\noindent Similarly, the spectrum of $\mathcal{L}^{+}$ consists of
curves: \begin{equation}
x=-\frac{y^{2}}{(c^{*})^{2}}-(1-\alpha)(1-k+\alpha k),\label{eq:3.13}\end{equation}
 or \begin{equation}
x=-\frac{y^{2}}{(c^{*})^{2}}+\alpha\label{eq:3.14}\end{equation}
 in the complex plane. Consequently, we have

\[
\max Re\,\sigma_{e}(\mathcal{L})=\max\{-\alpha,\,-1,\,-(1-\alpha)(1-k+\alpha k),\,\alpha\}=\alpha>0.\]

Hence, by \cite{07-Henry,20-Sandstede}, the steady state $U^{*}(\xi)$
of \eqref{eq:3.01} is unstable in $C_{0}$. \vspace{5pt}

\end{proof}
\noindent Theorem \ref{thm:9} says that the $C_{0}$ norm for the
initial perturbation is too large to stabilize the traveling wave
solutions. To have further control on the exponential rates of functions
in $C_{0}$, we introduce the following weighted Banach space. Let
$\sigma_{1}$, $\sigma_{2}$ be two non-negative numbers and the space
$C_{\sigma_{1},\sigma_{2}}$ be defined as: \[
C_{\sigma_{1},\sigma_{2}}=\{U(\xi)\mid U(\xi)(e^{\sigma_{1}\xi}+e^{-\sigma_{2}\xi})\in C_{0}\},\]
 with the norm \[
\parallel U\parallel_{C_{\sigma_{1},\sigma_{2}}}=\sup_{\xi\in\mathbb{R}}\parallel U(\xi)(e^{\sigma_{1}\xi}+e^{-\sigma_{2}\xi})\parallel.\]

\noindent We can also similarly define the Banach space $C_{\sigma_{1},\sigma_{2}}^{(2)}$
as:\[
C_{\sigma_{1},\sigma_{2}}^{(2)}=\{U|\,\, U(\cdot),U'(\cdot),U''(\cdot)\in C_{\sigma_{1},\sigma_{;2}}\}\]
 with norm\[
\left\Vert U\right\Vert _{C_{\sigma_{1},\sigma_{2}}^{(2)}}=\sup_{\xi\in\mathbb{R}}\Sigma_{i=0}^{2}\left\Vert (e^{\sigma_{1}\xi}+e^{-\sigma_{2}\xi})\frac{d^{i}U(\xi)}{d\xi}\right\Vert .\]
 It can be easily verified that $C_{\sigma_{1},\sigma_{2}}$ and $C_{\sigma_{1},\sigma_{2}}^{(2)}$
are both Banach spaces.

We now restrict the initial conditions and the operator $\mathcal{L}$
to the newly defined Banach space $C_{\sigma_{1},\sigma_{2}}$ with
$\sigma_{1}\geq0$, $\sigma_{2}\geq0$ and $\sigma_{1}^{2}+\sigma_{2}^{2}\neq0$.
To relate the operator $\mathcal{L}$ in $C_{0}$ to an equivalent
operator in $C_{\sigma_{1},\sigma_{2}}$, we introduce the following
translation operator $\mathcal{T}:\, C_{\sigma_{1},\sigma_{2}}\rightarrow C_{0}$
defined as \[
\mathcal{T}V:=(e^{\sigma_{1}\xi}+e^{-\sigma_{2}\xi})V.\]
 The operator $T$ is thus linear, bounded and has a bounded inverse
$\mathcal{T}^{-1}:C_{0}\rightarrow C_{\sigma_{1},\sigma_{2}}$ with
$\mathcal{T}^{-1}V=(e^{\sigma_{1}\xi}+e^{-\sigma_{2}\xi})^{-1}V$.

Consider operator \begin{equation}
\tilde{\mathcal{L}}V=\mathcal{T}\mathcal{L}\mathcal{T}^{-1}V.\label{eq:3.15}\end{equation}
 One can easily see that $\tilde{\mathcal{L}}$ is a linear operator
with domain $C^{(2)}(\mathbb{R})\times C^{(2)}(\mathbb{R})$. By relation
(\ref{eq:3.15}), considering $\mathcal{L}$ in $C_{\sigma_{1},\sigma_{2}}$
is equivalent to considering $\tilde{\mathcal{L}}$ in $C_{0}$. The
operator $\tilde{\mathcal{L}}$ can be explicitly written as \begin{equation}
\tilde{\mathcal{L}}V=V_{\xi\xi}-(2g_{1}+c^{*})V_{\xi}+(2g_{1}^{2}-g_{2}+c^{*}g_{1}+\frac{\partial F}{\partial U}(U^{*}))V,\label{eq:3.16}\end{equation}
 where $(2g_{1}+c^{*})$ and $(2g_{1}^{2}-g_{2}+c^{*}g_{1})$ in \eqref{eq:3.16}
represent the matrices $(2g_{1}+c^{*})I$ and $(2g_{1}^{2}-g_{2}+c^{*}g_{1})I$
with \[
\begin{array}{lll}
g_{1}(\xi) & = & {\displaystyle \frac{\sigma_{1}e^{\sigma_{1}\xi}-\sigma_{2}e^{-\sigma_{2}\xi}}{e^{\sigma_{1}\xi}+e^{-\sigma_{2}\xi}}},\\
\\g_{2}(\xi) & = & {\displaystyle \frac{\sigma_{1}^{2}e^{\sigma_{1}\xi}+\sigma_{2}^{2}e^{-\sigma_{2}\xi}}{e^{\sigma_{1}\xi}+e^{-\sigma_{2}\xi}}},\end{array}\]
 and further properties \[
\begin{array}{ll}
\lim_{\xi\rightarrow\infty}g_{1}(\xi)=\sigma_{1},\quad & \lim_{\xi\rightarrow-\infty}g_{1}(\xi)=-\sigma_{2};\\
\\\lim_{\xi\rightarrow\infty}g_{2}(\xi)=\sigma_{1}^{2}, & \lim_{\xi\rightarrow-\infty}g_{2}(\xi)=\sigma_{2}^{2}.\end{array}\]
 For the convenience of later use let $M(\xi)\doteq2g_{1}^{2}-g_{2}+c^{*}g_{1}+\frac{\partial F}{\partial U}(U^{*})$.
We now locate the essential spectrum of the operator $\tilde{\mathcal{L}}$
in the space $C_{0}$.
\begin{lem}
\textbf{\label{lem:10}} Let $c^{*}>2\sqrt{\alpha}$ and if $\sigma_{1}$
and $\sigma_{2}$ satisfy \begin{equation}
0\leq\sigma_{1}<\frac{-c^{*}+\sqrt{c^{*2}+4\alpha}}{2},\quad\frac{c^{*}-\sqrt{c^{*2}-4\alpha}}{2}<\sigma_{2}<\frac{c^{*}+\sqrt{c^{*2}-4\alpha}}{2},\label{eq:3.17}\end{equation}
then the essential spectrum of the operator $\tilde{\mathcal{L}}$
in the space $C_{0}(\mathbb{R})\times C_{0}(\mathbb{R})$ is contained
in a closed sector in the left half complex plane with vertex on the
horizontal axis left to the origin. Outside this sector, there are
only eigenvalues and the resolvent of $\tilde{\mathcal{L}}$. \end{lem}
\begin{proof}
\noindent Same as in the proof of Theorem \ref{thm:9}, we first study
the operator $\tilde{\mathcal{L}}$ at infinities. We have \begin{equation}
\tilde{\mathcal{L}}^{+}V=V_{\xi\xi}-(2\sigma_{1}+c^{*})V_{\xi}+(\sigma_{1}^{2}+c^{*}\sigma_{1}+\frac{\partial F}{\partial U}(U^{*}))V,\label{eq:3.18}\end{equation}
 \begin{equation}
\tilde{\mathcal{L}}^{-}V=V_{\xi\xi}-(-2\sigma_{2}+c^{*})V_{\xi}+(\sigma_{2}^{2}-c^{*}\sigma_{2}+\frac{\partial F}{\partial U}(U^{*}))V,\label{eq:3.19}\end{equation}
 where $\sigma_{1}^{2}+c^{*}\sigma_{1}+\frac{\partial F}{\partial U}(U^{*})$
and $\sigma_{2}^{2}-c^{*}\sigma_{2}+\frac{\partial F}{\partial U}(U^{*})$
are, respectively, the constant matrices \begin{equation}
M^{+}=\left[\begin{array}{cc}
\sigma_{1}^{2}+c^{*}\sigma_{1}-\alpha & 0\\
\\1-k & \sigma_{1}^{2}+c^{*}\sigma_{1}-1\end{array}\right],\label{eq:3.20}\end{equation}
 and \begin{equation}
M^{-}=\left[\begin{array}{cc}
\sigma_{2}^{2}-c^{*}\sigma_{2}-(1-\alpha)(1-k+k\alpha) & 1-\alpha\\
\\0 & \sigma_{2}^{2}-c^{*}\sigma_{2}+\alpha\end{array}\right].\label{eq:3.21}\end{equation}
 The essential spectrum of the operator $\tilde{\mathcal{L}}^{+}$
consists of curves:\[
\det(-\zeta^{2}I-i(2\sigma_{1}+c^{*})I\zeta+M^{+})=0.\]

\noindent Therefore, \[
-\zeta^{2}-i(2\sigma_{1}+c^{*})\zeta+\sigma_{1}^{2}+c^{*}\sigma_{1}-\alpha-\lambda=0,\]
 or \[
-\zeta^{2}-i(2\sigma_{1}+c^{*})\zeta+\sigma_{1}^{2}+c^{*}\sigma_{1}+\epsilon_{1}-1-\lambda=0.\]
 Equivalently , one has

\noindent \[
x=-\frac{y^{2}}{(2\sigma_{1}+c^{*})^{2}}+\sigma_{1}^{2}+c^{*}\sigma_{1}-\alpha\]
 or

\[
x=-\frac{y^{2}}{(2\sigma_{1}+c^{*})^{2}}+\sigma_{1}^{2}+c^{*}\sigma_{1}+\epsilon_{1}-1.\]

\noindent Similarly, the essential spectrum of $\mathcal{L}^{-}$
consists of curves: \[
x=-\frac{y^{2}}{(-2\sigma_{2}+c^{*})^{2}}+\sigma_{2}^{2}-c^{*}\sigma_{2}-(1-\alpha)(1-k+k\alpha),\]
 or \[
x=-\frac{y^{2}}{(-2\sigma_{2}+c^{*})^{2}}+\sigma_{2}^{2}-c^{*}\sigma_{2}+\alpha.\]
 Consequently, we have

\begin{equation}
\begin{array}{lll}
\max\{\mbox{Re}\,\sigma_{e}(\tilde{\mathcal{L}})\} & = & \max\{\sigma_{1}^{2}+c^{*}\sigma_{1}-\alpha,\,\sigma_{1}^{2}+c^{*}\sigma_{1}-1,\\
\\ &  & \sigma_{2}^{2}-c^{*}\sigma_{2}{\displaystyle -(1-\alpha)(1-k+k\alpha),\,\sigma_{2}^{2}-c^{*}\sigma_{2}+\alpha\}}\end{array}\label{eq:3.22}\end{equation}
 A simple calculation shows the validity of the Lemma. \end{proof}
\begin{cor}
\textbf{\label{pro:-Corollary-11}} Assume that all the hypotheses
in Lemma \ref{lem:10} hold, then the essential spectrum of the operator
$\mathcal{L}$ in the space $C_{\sigma_{1},\sigma_{2}}$ lies in the
left half complex plane. \vspace{5pt}
\end{cor}
\begin{proof}
\noindent The conclusion follows immediately from Lemma \ref{lem:10}
and relation (\ref{eq:3.15}). \vspace{5pt}

\end{proof}
We now locate the eigenvalues of the operator $\mathcal{L}$ in the
space $C_{\sigma_{1},\sigma_{2}}$. The next lemma shows that $0$
is not an eigenvalue of the operator $\mathcal{L}$ in the space $C_{\sigma_{1},\sigma_{2}}$.
The idea of the proof is as follows. If there is an eigenfunction
associated with the eigenvalue zero of the operator $\mathcal{L}$
in the space $C_{\sigma_{1},\sigma_{2}}$, then we can show the derivative
of the traveling wave solution is larger than or equal to any multiple
of this function by examing their asymptotics at infinities. But the
boundedness of the derivative of the traveling wave solution makes
this impossible.
\begin{lem}
\textbf{\label{lem:12}} Assuming all the hypotheses in Lemma \ref{lem:10},
then $0$ is not an eigenvalue of the operator $\mathcal{L}$ in $C_{\sigma_{1},\sigma_{2}}$.
\vspace{5pt}
\end{lem}
\begin{proof}
\noindent Let $U^{*}(\xi)$ be a traveling wave solution of \eqref{eq:3.01}
found in Theorem \ref{thm:2.06}. By taking Gateaux derivative we
can verify that \begin{equation}
\mathcal{L}((U^{*}(\xi))')=0.\label{eq:3.23}\end{equation}

Suppose that there exists a nonzero function $\bar{V}\in C_{\sigma_{1},\sigma_{2}}$
satisfying the equation \begin{equation}
\mathcal{L}\bar{V}=0.\label{eq:3.24}\end{equation}
we then claim that the inequality $|r\bar{V}(\xi)|\leq(U^{*})'(\xi)$
is consequently true for all $r\in\mathbb{R}$ and $\xi\in\mathbb{R}$.
Writing $S=\{r\in\mathbb{R}|\,|r\bar{V}(\xi)|\leq(U^{*})'(\xi),\:\xi\in\mathbb{R}\}$,
we verify the following properties of $S$:

1. $S$ is non-empty since $0\in S$.

2. $S$ is closed. Let $r_{i}\in S$, $i=1,\,2,\,...$ and $r_{i}\rightarrow r$
as $i\rightarrow+\infty$, then we will have $|r_{i}\bar{V}(\xi)|\leq(U^{*})'$,
which implies that $|r\bar{V}(\xi)|\leq(U^{*})'(\xi)$, we therefore
have $r\in S$.

3. $S$ is open. Let $r\in S$, we will show that there exists a $\bar{\delta}>0$
such that $(r-\bar{\delta},\, r+\bar{\delta})\subset S$. We \textit{claim}
that $|r\bar{V}(\xi)|\leq(U^{*})'(\xi)$ implies $|r\bar{V}(\xi)|<(U^{*})'(\xi)$.
In fact, let $W(\xi)=(U^{*})'(\xi)-r\bar{V}(\xi)$ then $W(\xi)\geq0$
for $\xi\in\mathbb{R}$ and it satisfies the following equation:\begin{equation}
\left\{ \begin{array}{l}
w_{1}''-cw_{1}'+A_{11}w_{1}+A_{12}w_{2}=0,\\
\\w_{2}''-cw_{2}'+A_{21}w_{1}+A_{22}w_{2}=0,\\
\\(w_{1},w_{2})(-\infty)=(w_{1},w_{2})(+\infty)=0,\end{array}\right.\label{eq:3.25}\end{equation}
where $A_{ij},$ $i,j=1,2$ are the entries of the Jacobian $\frac{\partial F}{\partial U}(U^{*})$.
The Maximum Principle implies that $W(\xi)=(w_{1}(\xi),w_{2}(\xi))^{T}>0$
for $\xi\in\mathbb{R}$. This shows that $(U^{*})'(\xi)-r\bar{V}(\xi)>0,\:\xi\in\mathbb{R}$.
Similarly we can show that $(U^{*})'(\xi)+r\bar{V}(\xi)>0$ for $\xi\in\mathbb{R}$.
The \textit{claim} then follows.

We next show that the \textit{claim} further implies $|\bar{r}\bar{V}(\xi)|<(U^{*})'(\xi)$
as long as $\bar{r}$ is sufficiently close to $r$. According to
condition \eqref{eq:3.17} and the assumption that $\bar{V}\in C_{\sigma_{1},\sigma_{2}}$,
for any fixed $\tilde{r}\in\mathbb{R}$, there exists a $N>0$ sufficiently
large such that $(e^{\sigma_{1}\xi}+e^{-\sigma_{2}\xi})[(U^{*})'(\xi)-\tilde{r}\bar{V}(\xi)]>0$
for all $\xi\in(-\infty,-N]$, this implies that $(U^{*})'(\xi)>\tilde{r}\bar{V}(\xi)$
also holds there. Furthermore, due to the \textit{claim} and the boundedness
of the functions $(U^{*})'$ and $\bar{V}$, we can find a $\bar{\delta}>0$
such that for any $\bar{r}\in(-\bar{\delta}+r,\,\bar{\delta}+r)$,
$(U^{*})'(\xi)>\bar{r}\bar{V}(\xi)$ on the finite interval $[-N,\, N]$.
Now we fix $\tilde{r}=\bar{r}$ and show $[(U^{*})'(\xi)-\bar{r}\bar{V}(\xi)]>0$
for $\xi\in[N,+\infty)$.

Noting that the entries of the diagonal of the matrix $\frac{\partial F}{\partial U}(U^{*}(+\infty))$
are both negative, we can choose a vector $P_{+}>0$ such that (by
increasing $N$ if necessary) $\frac{\partial F}{\partial U}(U^{*}(\xi))P_{+}<0$
for $\xi\in[N,+\infty)$.

We have to deal with two cases:

Case A. If we already have $(U^{*})'(\xi)-\bar{r}\bar{V}(\xi)\geq0$
for $\xi\geq N$, then the Maximum Principle implies that $(U^{*})'(\xi)-\bar{r}\bar{V}(\xi)>0$
on $[N,+\infty)$. Analogly $(U^{*})'(\xi)+\bar{r}\bar{V}(\xi)>0$
is also true for $\xi\in\mathbb{R}$. Consequently, $S$ is open.

Case B. If there is a point in the interval $(-\infty,-N)$ such that
one of the components of vector $(U^{*})'(\xi)-\bar{r}\bar{V}(\xi)$
takes negative local minimum on this point, we consider function $\tilde{W}(\xi)=(U^{*})'(\xi)-\bar{r}\bar{V}(\xi)+\tau P_{+}$.
The asymptotic rates of $(U^{*})'$ and $\bar{V}$ at $+\infty$ imply
that there is a sufficiently large $\tau>0$ such that $\tilde{W}=(U^{*})'(\xi)-\bar{r}\bar{V}(\xi)+\tau P_{+}\geq0$
for $\xi\in[N,+\infty)$. We further assume that one of the components
of $\tilde{W}(\xi)$, say $\tilde{w}_{1}$ for example, takes minimum
at a finite point $\xi_{2}\,\mbox{in}\,[N,+\infty)$. It is not hard
to verify that there is a $\tau=\tau_{2}>0$ such that the corresponding
$\tilde{W}(\xi)$ staisfying $\tilde{w}_{1}(\xi_{2})=0$ and $\tilde{W}(\xi)\geq0$
for $\xi\in(-\infty,-N)$. For such $\tau_{2}$ on the one hand, we
have\begin{equation}
\begin{array}{lll}
\mathcal{L}\tilde{W} & = & \tilde{W}_{\xi\xi}-c^{*}\tilde{W}_{\xi}+\frac{\partial F}{\partial U}(U^{*})\tilde{W}\\
\\ & = & \tau_{2}\frac{\partial F}{\partial U}(U^{*})P_{+}<0,\end{array}\label{eq:3.26}\end{equation}
on the other hand at $\xi=\xi_{2}$, the first component to the left
hand side of \eqref{eq:3.26} is larger than or equals zero. We then
have a contradiction, and consequently $(U^{*})'(\xi)-\bar{r}\bar{V}(\xi)\geq0$
for $\xi\in[N,+\infty)$. We are again in the situation descibed by
case A. By a similar argument, we can show $(U^{*})'(\xi)+\bar{r}\bar{V}(\xi)\geq0$
for  $\xi\in[N,+\infty)$.

In summary, for any $\bar{r}\in(-\bar{\delta}+r,\,\bar{\delta}+r)$,
we have $|\bar{r}\bar{V}(\xi)|<(U^{*})'$, i.e., $S$ is open. 

Now the set $S$ is a non-empty, open and closed subset of $\mathbb{R}$,
then $S\equiv\mathbb{R}$, but this is impossible since $(U^{*})'$
is bounded. Therefore the equation \eqref{eq:3.24} does not have
non-zero solution in $C_{\sigma_{1},\sigma_{2}}$.
\end{proof}
The next lemma shows that there is no eigenvalue of $\mathcal{L}$
in $C_{\sigma_{1},\sigma_{2}}$ with positive or zero real part.
\begin{lem}
\label{lem:13} Let $C_{0}^{\mathbb{C}}$ be the complexified space
of $C_{0}(\mathbb{R})$ and $\lambda$ be an eigenvalue of the operator
$\tilde{\mathcal{L}}$ with $\underline{U}\in C_{0}^{\mathbb{C}}$
as its eigenfunction, then $\mbox{Re}\:\lambda<0$. \end{lem}
\begin{proof}
Let $\lambda=\lambda_{1}+\lambda_{2}i$ and $\underline{U}(\xi)=U^{1}(\xi)+iU^{2}(\xi)$
for $\xi\in\mathbb{R}$, where $\lambda_{1},\,\lambda_{2}\in\mathbb{R}$
with $\lambda_{1}^{2}+\lambda_{2}^{2}\neq0$ and $U^{i}(\xi)\in C_{0}(\mathbb{R})$. 

Consider the Cauchy problem (\cite{03-BatesChen,28-XuZhao}),

\begin{equation}
V_{t}=\tilde{\mathcal{L}}V-\lambda_{1}V,\quad V(\xi,0)=U^{1}(\xi).\label{eq:3.27}\end{equation}
 It is easy to verify that $V(\xi,t)=U^{1}(\xi)\cos(\lambda_{2}t)-U^{2}(\xi)\sin(\lambda_{2}t)$
is bounded and solves \eqref{eq:3.27} for $\xi\in\mathbb{R}$ and
$t\geq0$. We may assume that at least one of the components of $V$
takes positive value at some point of $\xi$ and $t$, otherwise we
can consider $-V$. Since the vector $\mathcal{T}(U^{*})'(\xi)\rightarrow+\infty$
as $\xi\rightarrow-\infty$, similar to the proof of Lemma \ref{lem:12},
we can choose a large $\xi_{0}>0$ such that

\begin{equation}
M(\xi)P^{+}<0\label{eq:3.28}\end{equation}
for $\xi\geq\xi_{0}$, and at the same time

\begin{equation}
V(\xi,t)<\mathcal{T}(U^{*})'(\xi)\label{eq:3.29}\end{equation}
for $\xi\leq-\xi_{0}$ and $t\geq0$, and $M(\xi)=2g_{1}^{2}-g_{2}+c^{*}g_{1}+\frac{\partial F}{\partial U}(U^{*})$
in \eqref{eq:3.28}. Moreover, we can choose a small $\bar{\epsilon}>0$
such that

\begin{equation}
(\bar{\epsilon}^{2}I-\bar{\epsilon}(2g_{1}+c^{*})I+M(\xi))P^{+}<0.\label{eq:3.30}\end{equation}
 for $\xi\geq\xi_{0}$. We can choose an even larger $\xi_{0}$ if
necessary such that at least one of the positive values of $V$ occurs
inside the interval $(-\xi_{0},\xi_{0})$. The positivity of $\mathcal{T}(U^{*})'(\xi)$
for $\xi\in\mathbb{R}$ implies that there is a $r>0$ such that

\begin{equation}
V(\xi,t)\leq r\mathcal{T}(U^{*})'(\xi)\quad\mbox{for}\:\xi\in[-\xi_{0},\xi_{0}],\, t\geq0.\label{eq:3.31}\end{equation}
 Suitably adjusting $r$ such that the equality in \eqref{eq:3.31}
holds on at least one component at a point $(\xi_{1},t_{1})$ with
$\xi_{1}\in[-\xi_{0},\xi_{0}]$ and $t_{1}\geq0$. Now let $\lambda_{1}\geq0$,
we have

\textit{Claim:} $V(\xi,t)\leq r\mathcal{T}(U^{*})'$ for all $\xi\in\mathbb{R}$,
$t\geq0$.

The claim is already true for $\xi\leq\xi_{0}$ by the discussion
above. Suppose the assertion is not true for $\xi>\xi_{0}$, then
we can find a $\bar{\xi}>\xi_{0}$, and a $t_{1}\geq0$ such that
$V(\bar{\xi,}t_{1})>r\mathcal{T}(U^{*})'(\bar{\xi})$. Since $Q^{+}(\xi)\doteq e^{\bar{\epsilon}\xi}P^{+}\rightarrow+\infty$
as $\xi\rightarrow+\infty$, there is a $\tilde{r}>0$ such that $V(\xi,t)\leq r\mathcal{T}(U^{*})'(\xi)+\tilde{r}Q^{+}(\xi)$
for all $\xi\geq\xi_{0}$ and $t\geq0$, and for at least one $j$
($j=1,2$ ) and one $\xi_{2}\geq\xi_{0}$ and one $t_{2}>0$, we have
the equality for the $j-$th component:

\[
V_{j}(\xi_{2},t_{2})=r\mathcal{T}(U_{j}^{*})'(\xi_{2})+\tilde{r}Q_{j}^{+}(\xi_{2}).\]
 Let $Y(\xi,t)=r\mathcal{T}(U^{*})'(\xi)+\tilde{r}Q^{+}(\xi)-V(\xi,t)$\textbf{,}
then $Y_{j}$ has the following properties: $Y_{j}(\xi_{2},t_{2})=0,$
$Y_{j}(\xi_{0},t)>0$, $Y_{j}(\xi,t)\geq0$ for all $\xi\geq\xi_{0}$,
$t\geq0$. Hence $Y_{j,t}(\xi_{2},t_{2})\leq0$, $Y_{j,\xi}(\xi_{2},t_{2})=0$
and $Y_{j,\xi\xi}(\xi_{2},t_{2})\geq0$, and $Y_{j}(\xi,t)$ satisfies

\begin{equation}
\begin{array}{lll}
Y_{j,t} & = & -V_{j,t}\\
 & = & -(\tilde{\mathcal{L}}V-\lambda_{1}V)_{j}\\
 & > & (-\tilde{\mathcal{L}}V+\lambda_{1}V+\tilde{\mathcal{L}}r\mathcal{T}(U^{*})'+\tilde{\mathcal{L}}\tilde{r}Q^{+}-\lambda_{1}(r\mathcal{T}(U^{*})'+\tilde{r}Q^{+}))_{j}\\
 & = & (\tilde{\mathcal{L}}Y-\lambda_{1}Y)_{j}\\
 & = & Y_{j,\xi\xi}-(2g_{1}+c^{*})Y_{j,\xi}+M_{j,1}Y_{1}+M_{j,2}Y_{2}-\lambda_{1}Y_{j}.\end{array}\label{eq:3.32}\end{equation}
owing to \eqref{eq:3.30}, $\tilde{\mathcal{L}}r\mathcal{T}(U^{*})'=0$
and $r\mathcal{T}(U^{*})'+\tilde{r}Q^{+})\geq0$. However, we find
that at $(\xi_{2},t_{2})$ the left hand side of the inequality \eqref{eq:3.32}
is smaller than or equal to zero, while the right hand side is larger
than zero. We then have a contradiction. Therefore the claim is proved.

We next show that the above is also true even for $r=0$. In fact,
let $S=\{r\geq0|V(\xi,t)\leq r\mathcal{T}(U^{*})'(\xi)\;\mbox{for}\;\xi\in\mathbb{R},\, t\geq0\}$,
by the claim $S\neq\varnothing$, then $r_{0}\doteq\min S\geq0$.
We need to show $r_{0}=0$. Suppose on the contrary that we have $r_{0}>0$
, then \begin{equation}
V(\xi,t)\leq r_{0}\mathcal{T}(U^{*})'(\xi)\;\mbox{for}\;\xi\in\mathbb{R},\, t\geq0.\label{eq:3.33}\end{equation}
 Since $\mathcal{T}(U^{*})'(\xi)\rightarrow+\infty$ as $\xi\rightarrow-\infty$,
the strict inequality holds at $\xi$ near $-\infty$. We first assume
that an equality occurs in \eqref{eq:3.33} at one of the components
(say, the $i$-th component) at a point $(\tilde{\xi},\tilde{t})$
with $\tilde{\xi}\in\mathbb{R}$ and $\tilde{t}\geq0$. Moreover from

\[
\begin{array}{lll}
Y_{i,t} & \geq & (\tilde{\mathcal{L}}Y-\lambda_{1}Y)_{i}\\
 & = & Y_{i,\xi\xi}-(2g_{1}+c^{*})Y_{i,\xi}+M_{i,1}(\xi)Y_{1}+M_{i,2}(\xi)Y_{2}-\lambda_{1}Y_{i}\end{array}\]
with $Y(\xi,t)=r(U^{*})'(\xi)-V(\xi,t)$ and $(M_{i,1},M_{i,2})$
the $i$-th row of the matrix function $M(\xi)$. The positivity Lemma
(\cite{23-Volpert}) implies that $Y_{i}(\xi,t)>0$ for $\xi\in\mathbb{R}$
and $t>\tilde{t}$, further by the periodicity of $V$, we have that
$Y_{i}(\xi,t)>0$ for $\xi\in\mathbb{R}$ and $t>0$. This contradicts
with the existence of $\tilde{\xi}$ and shows that

\[
V(\xi,t)<r_{0}\mathcal{T}(U^{*})'(\xi)\;\mbox{for}\;\xi\in\mathbb{R},\: t>0.\]

Similar to the justification of the \textit{claim} in the proof of
Lemma \ref{lem:12}, we can show there exists a sufficiently small
$\delta$ with $r_{0}\gg\delta>0$, such that

\[
V(\xi,t)<(r_{0}-\delta)\mathcal{T}(U^{*})'(\xi)\;\mbox{for}\;\xi\in\mathbb{R},\: t>0.\]
It then follows that $r_{0}-\delta\in S$ thus contradicts with the
definition of $r_{0}$. Then we conclude that $r_{0}=0$ but this
again is impossible since by assumption a component of $V$ is positive
at some point $\xi$ and $t$. Hence we have proved the Lemma. \end{proof}
\begin{lem}
\label{lem:14}Under the hypotheses of Lemma \ref{lem:10}, the operator
$\tilde{\mathcal{L}}$ in $C_{0}$ is sectorial. \end{lem}
\begin{proof}
\noindent See \cite{07-Henry}.\end{proof}
\begin{thm}
\noindent \label{thm:15}Under the hypotheses of \ref{lem:10}, the
operator $\mathcal{L}$ generates an analytical semigroup in $C_{\sigma_{1},\sigma_{2}}$.\end{thm}
\begin{proof}
The conclusion follows from Lemma \ref{lem:14} and \cite{07-Henry}. 
\end{proof}
We now state the stability theorem,
\begin{thm}
\label{thm:16}Assume the hypotheses of Lemma \ref{lem:10}. The traveling
wave solution $U^{*}$ of \eqref{eq:1.02}, with wave speed $c^{*}>2\sqrt{\alpha}$
, is asymptotically stable according to norm $||\cdot||:=||\cdot||_{C_{\sigma_{1},\sigma_{2}}}$.
That is, if the initial condition $U(\xi,0)=\bar{U}(\xi)\in C$ with
$(\bar{U}(\xi)-U^{*}(\xi))\in C_{\sigma_{1},\sigma_{2}}$ and $||\bar{U}-U^{*}||$
sufficiently small, then the solution $U(\xi,t)$ to \eqref{eq:3.01}
exists uniquely for all $t>0$ and satisfies \begin{equation}
||U(\xi,t)-U^{*}(x+ct)||\leq Me^{-bt},\label{eq:3.34}\end{equation}
where the constants $M>0,\, b>0$ are independent of $t$ and $\bar{U}$. \end{thm}
\begin{proof}
Corollary \ref{pro:-Corollary-11}, Lemma \ref{lem:10}, Lemma \ref{lem:12}
and \ref{lem:13} show that the spectrum of the operator $\mathcal{L}$
in the space $C_{\sigma_{1},\sigma_{2}}$ is contained in an angular
region in the left complex plane. Thus, the trivial solution of system
\eqref{eq:3.01} is linearly asymptotically stable. Theorem \ref{thm:15}
further implies the local nonlinear stability of the traveling wave
solutions (\cite{02-AlexanderGardnerJones,17-Kapitula,20-Sandstede}). 
\end{proof}
\textbf{\textit{Acknowledgement: }}The authors thank Professor Anthony
W. Leung of University of Cincinnati for the inspiration of the lower
solution derived in Lemma \ref{lem:2.04}.

\end{document}